\documentclass[12pt]{elsarticle}
\usepackage{amsmath,amssymb,amsfonts, hyperref}

\newtheorem{theorem}{Theorem}
\newtheorem{lemma}{Lemma}

\newtheorem{claim}{Claim}
\newtheorem{conjecture}{Conjecture}

\newcommand{\B}{{\mathrm{B}}}

\begin{document}

\begin{center}
{\Large
On monotonicity of Ramanujan function for binomial random variables}\\

\vspace{0.2cm}

{\large Daniil Dmitriev\footnote{\'{E}cole Polytechnique F\'{e}d\'{e}rale de Lausanne, laboratory of machine learning and optimization\\ daniil.dmitriev@epfl.ch},
 Maksim Zhukovskii\footnote{Sobolev Institute of Mathematics, Pevtsova, 13, Omsk, 644043, Russia; Moscow Institute of Physics and Technology, laboratory of advanced combinatorics and network applications, 9 Institutskiy per., Dolgoprudny, Moscow Region, 141701, Russian Federation; The Russian Presidential Academy of National Economy and Public Administration, Prospekt Vernadskogo, 82, Moscow, Russian Federation, 119571; Adyghe State University, Caucasus mathematical center, Ulitsa Pervomayskaya, 208, Maykop, Respublika Adygeya, 385000\\ zhukmax@gmail.com}
 }
\vspace{0.5cm}

Abstract

\end{center}

For a binomial random variable $\xi$ with parameters $n$ and $b/n$, it is well known that the median equals $b$ when $b$ is an integer. In 1968, Jogdeo and Samuels studied the behaviour of the relative difference between ${\sf P}(\xi=b)$ and $1/2-{\sf P}(\xi<b)$. They proved its monotonicity in $n$ and posed a question about its monotonicity in $b$. This question is motivated by the solved problem proposed by Ramanujan in 1911 on the monotonicity of the same quantity but for a Poisson random variable with an integer parameter $b$. In the paper, we answer this question and introduce a simple way to analyse the monotonicity of similar functions.

\vspace{0.5cm}
\vspace{0.2cm}
\section{Introduction}

Given a non-negative integer random variable $\xi$, we call {\it the median of $\xi$}
$$
 \mu(\xi):=\min\left\{m\in\mathbb{Z}_+:\,{\sf P}(\xi\leq m)\geq\frac{1}{2}\right\}.
$$

Consider a Poisson random variable $\eta_b$ with a positive integer parameter $b$. From the general result of Choi~\cite{Choi} it follows that $\mu(\eta_b) = b$. 
The more delicate question concerns the behaviour of the distribution function of $\eta_b$ around $b$ and was initially posed by S. Ramanujan~\cite{Ram}, who conjectured that
$$
y_b:=\frac{\frac{1}{2}-{\sf P}(\eta_b<b)}{{\sf P}(\eta_b=b)}\in\left(\frac{1}{3},\frac{1}{2}\right)
$$
and $y_b$ decreases. This was proven independently by G. Szeg\"{o} in 1928~\cite{Szego} and G.N. Watson in 1929~\cite{Watson} and also implies that $\mu(\eta_b)=b$. Since then, the behaviour of $y_b$ was widely studied. Below, we give a very brief history of this study.

In 1913, in his letter to Hardy, Ramanujan posed an initial question: he conjectured that
$$
y_b=\frac{1}{3}+\frac{4}{135(b+\alpha_b)},
$$
where $8/45\geq\alpha_b\geq 2/21$. This conjecture was proved by Flajolet et al.~\cite{Flajolet} in 1995. In 2003, S. E. Alm~\cite{Alm} proved that $\alpha_b$ decreases. In 2004, this result was strengthened by H. Alzer~\cite{Alzer}:
$$
 y_b=\frac{1}{3}+\frac{4}{135b}-\frac{8}{2835(b^2+\beta_b)},
$$
where $-\frac{1}{3}<\beta_b\leq-1+\frac{4}{\sqrt{21(368-135e)}}$, and the bounds are sharp.\\

Similar properties are of interest for a binomial random variable $\xi_{b,n}$ with parameters $n$ and $b/n$, where $b\leq n$ are positive integers. It is known~\cite{Hamza} that $\left|\mu(\xi_{b,n})-b\right|\leq\ln 2$, and so, $\mu(\xi_{b,n})=b$. In 1968, K. Jogdeo and S. M. Samuels~\cite{Samuels_med} studied the generalization of the first Ramanujan conjecture for the binomial random variables. They considered the behaviour of
\begin{equation}
\label{z_bn_def}
z_{b,n}:=\frac{\frac{1}{2}-{\sf P}(\xi_{b,n}<b)}{{\sf P}(\xi_{b,n}=b)},
\end{equation}
and proved the following result.
\begin{theorem}[K. Jogdeo, S. M. Samuels, 1968]
For every $b$, $z_{b,n}$ decreases for $n\geq 2b$ and $z_{b,n}\to y_n$ as $n\to\infty$. Moreover, for all $n>2b$, $1/3<z_{b,n}<1/2$. For all $b<n<2b$, $1/2<z_{b,n}<2/3$. Finally, $z_{b,2b}=1/2=z_{b,b}$.
\label{Sam}
\end{theorem}
In the paper, they also mentioned that, obviously, $z_{b+1,n}<z_{b,n}$ for all large enough $n$, but they were unable to make this more precise.

We solve the problem proposed by Jogdeo and Samuels on a monotonicity of $z_{b,n}$ in $b$. Our main result is the following.
\begin{theorem}
Let $\varepsilon>0$. There exists $n_0$ such that, for all $n\geq n_0$, 
\begin{enumerate}
\item if $n-(1+\varepsilon)\frac{2}{3\sqrt{5}}\sqrt{n}>b>(1+\varepsilon)\frac{2}{3\sqrt{5}}\sqrt{n}$, then $z_{b+1,n}>z_{b,n}$,
\item if either $b>n-(1-\varepsilon)\frac{2}{3\sqrt{5}}\sqrt{n}$ or $b<(1-\varepsilon)\frac{2}{3\sqrt{5}}\sqrt{n}$, then $z_{b+1,n}<z_{b,n}$.
\end{enumerate}
\label{zs}
\end{theorem}

We are also interested in a monotonicity of $p_{b,n}:={\sf P}(\xi_{b,n}<b)$ in $b$. From the result of  Szeg\"{o} and Watson, it immediately follows that ${\sf P}(\eta_b<b)$ increases (or, in other words, the difference between $1/2$ and the probability that $\eta_b$ is less than the median decreases). So, for $n$ large enough, $p_{b+1,n}>p_{b,n}$ as well. It is easy to see that for $n=b+1$, $0=p_{b+1,n}<p_{b,n}$, in contrast. For small values of $b$, it can be verified that the same inequality holds even for $b+1\leq n\leq 3b+1$.  We prove that for all values of $b$ the monotonicity of $p_{b,n}$ in $b$ changes when $n = 3b + 2$.

\begin{theorem}
The following properties hold.
\begin{itemize}
\item If $n\geq 3b+2$, then $p_{b+1,n}>p_{b,n}$.

\item If $ n\leq 3b+1$, then $p_{b+1,n}<p_{b,n}$.
\end{itemize}
\label{ps}
\end{theorem}
In Appendix D, for illustration, we give plots for $p_{b,n}$ and $z_{b,n}$.

The rest of the paper is organized in the following way. In Section~\ref{tools}, we describe the main tools. In Section~\ref{th2}, we prove Theorem~\ref{ps}. In Section~\ref{th1}, we prove Theorem~\ref{zs}. Section~\ref{further} is devoted to another motivation for our result, a certain inequality of small deviations closely related to Samuel's conjecture. 

\section{Main tools}
\label{tools}

The main ingredient of our proofs is avoiding summations in the definitions of $p_{b,n}$ and $z_{b,n}$ by replacing the factorials with gamma functions. It is known that it works well for Poisson distributions (see~\cite{Watson}). We show that it also helps to eliminate summations in the case of binomial distributions (see Claim~\ref{pg_expression}). Having this, we show that (see Section~\ref{2.1}) $g_{b,n}:=\int_{\Delta_{b,n}}g(z)dz$, where $g(z)=(1-z)^{b-1}z^{n-b}$ and $\Delta_{b,n}:=[1-(b+1)/n,1-b/n]$, gives the major contribution to $p_{b+1,n}-p_{b,n}$ (a similar observation for $z_{b+1,n}-z_{b,n}$ is also obtained, see Equations (\ref{p_difference_g}) and (\ref{g_useful}) of Claim~\ref{pg_expression}).

In Section~\ref{function_g},  we give very tight lower and upper bounds on $g$ based on Taylor expansion of $g$. It turns out that for our goals, it is sufficient to consider the first four terms in the Taylor polynomial.  Claim~\ref{pg_expression}, mentioned bounds in Section~\ref{function_g} and very careful analysis of obtained polynomials are sufficient for proving Theorem~\ref{ps}. This technique also works for proving Theorem~\ref{zs} when $c\sqrt{n}<b<n/2$ for an appropriate choice of $c>0$. For $b\leq c\sqrt{n}$, we exploit an asymptotical expansion of $z_{b,n}$ that is obtained in Section~\ref{g_asymptotics}. The case $b\geq n/2$ follows immediately from the observation that $z_{b,n}=1-z_{n-b,n}$.


All missing technical proofs (for Claims~\ref{pg_expression},~\ref{CL_below_0.5}~and~\ref{CL_z_estimation}) are provided in Appendix A.

\subsection{An integral expressions for $p_{b,n}$ and $z_{b, n}$}
\label{2.1}

The values $p_{b,n}$ by definition are equal to $\sum_{i=0}^{b-1}{n\choose i}(b/n)^i(1-b/n)^{n-i}$. In order to analyze this quantity we move from the discrete summation to the integral of the function $g(z)$. The same trick can be done for the value $z_{b,n}$, as stated in this result:

\begin{claim} 
\label{pg_expression}
For $p_{b,n}$ and consecutive differences $p_{b+1,n}-p_{b,n}$ the following equalities hold:
\begin{equation}
\label{from_p_to_beta}
p_{b,n} = \frac{\int_0^{1-b/n}(1-z)^{b-1}z^{n-b}dz}{\int_0^1(1-z)^{b-1}z^{n-b}dz},
\end{equation}
\begin{equation}
\label{p_difference_g}
p_{b+1,n}-p_{b,n}=\frac{\left(\frac{b+1}{n}\right)^b\left(1-\frac{b+1}{n}\right)^{n-b}-b\int_{1-(b+1)/n}^{1-b/n}(1-z)^{b-1} z^{n-b}dz}{b\int_0^1(1-z)^{b-1}z^{n-b}dz}.
\end{equation}
For $z_{b,n}$ and consecutive differences $z_{b+1,n}-z_{b,n}$ the following equalities hold:
\begin{equation}
\label{zs_integral}
z_{b,n} =\frac{\frac{1}{2}b\left(\left[\int_{1-b/n}^1-\int_{0}^{1-b/n}\right]g(z)dz\right)}{(b/n)^b(1-b/n)^{n-b}},
\end{equation}
$$
z_{b+1,n}-z_{b,n}=\frac{n^{n}}{2(n-b)(b+1)^b(n-b-1)^{n-b-1}}\times
\biggl(-2\left(\frac{b+1}{n}\right)^b\left(\frac{n-b-1}{n}\right)^{n-b}+
$$
$$
b\left[\int_{1-(b+1)/n}^1-\int_{0}^{1-(b+1)/n}\right]g(z)dz-
$$
\begin{equation}
\label{g_useful}
b\left(1+\frac{1}{b}\right)^b\left(1-\frac{1}{n-b}\right)^{n-b-1}\left[\int_{1-b/n}^1-\int_{0}^{1-b/n}\right]g(z)dz\biggr).
\end{equation}
\end{claim}

{\it Proof: see Appendix A.}

\subsection{Bounds on $g$ inside $\Delta_{b,n}$}
\label{function_g}

An inductive proof of the following observation is straightforward.

\begin{claim}
Let $\ell\in\{1,\ldots,\min\{b-1,n-b\}\}$. Then 
$$
\frac{\partial^{\ell} g}{\partial z^{\ell}}=(1-z)^{b-1-\ell}z^{n-b-\ell}\sum_{i=0}^{\ell}{\ell\choose i}(-1)^{\ell-i}z^{\ell-i}\frac{(n-1-i)!}{(n-1-\ell)!}\frac{(n-b)!}{(n-b-i)!}.
$$
\label{CL2}
\end{claim}

Using Claim~\ref{CL2} and the Lagrange form of the remainder of the Taylor polynomial, we obtain a lower and an upper bounds for $g(z)$ on $\Delta_{b,n}$.

Let, for $\ell\in\mathbb{Z}_+$, 
$
g_{\ell}(z)=\frac{1}{\ell!}\left(z-1+\frac{b+1}{n}\right)^{\ell}\left[\frac{\partial^\ell g}{\partial z^{\ell}}\left(1-\frac{b+1}{n}\right)\right]
$
be the $\ell$-th term in the Taylor expansion of $g$, let $g_{\leq 3}(z)=\sum_{\ell=0}^3 g_{\ell}(z)$, and 
$
d_4^-(z)=\frac{\partial^4 g}{\partial z^{4}}\left(1-\frac{b+1}{n}\right), d_4^+(z)=\frac{\partial^4 g}{\partial z^{4}}\left(1-\frac{b}{n}\right).
$

\begin{claim}
For $5\leq b\leq n/2$  and all $z\in\Delta_{b,n}$,
$$
g_{\leq 3}(z)+\frac{1}{24} \left(z-1+\frac{b+1}{n}\right)^4d_4^-(z)\leq g(z)\leq g_{\leq 3}(z)+ \frac{1}{24} \left(z-1+\frac{b+1}{n}\right)^4d_4^+(z).
$$
\label{CL_below_0.5}
\end{claim}
{\it Proof.} Immediately follows from the fact that $\partial^4 g/\partial z^4$ increases on $\Delta_{b,n}$. The latter is proven in Appendix A.

\subsection{Behaviour of $z_{b,n}$ for $b=O(\sqrt{n})$}
\label{g_asymptotics}

It is well known (see, e.g.,~\cite{Cheng}) that $y_b=\frac{1}{3}+\frac{4}{135b}-\frac{8}{2835b^2}+O\left(\frac{1}{b^3}\right).$ From this and Stirling's approximation (see, e.g.,~\cite{Feller}) $b!=\sqrt{2\pi b}b^b\exp\left(-b+\frac{1}{12b}+O(1/b^3)\right)$, it follows that, for a Poisson random variable $\eta_b$,
\begin{equation}
 {\sf P}(\eta_b<b)=\frac{1}{2}-\frac{1}{3\sqrt{2\pi b}}-\frac{1}{540\sqrt{2\pi}b\sqrt{b}}+O\left(\frac{1}{b^2\sqrt{b}}\right).
\label{p_b_expantion}
\end{equation}

For a non-negative integer $\xi$, we denote $\xi^{(s)}:=\xi(\xi-1)\ldots(\xi-s+1)$.

\begin{lemma}
For every $s\in\mathbb{N}$,
$$
{\sf E}\left[\eta_b^{(s)}I(\eta_b<b)\right]=
b^s\left({\sf P}(\eta_b<b)-\frac{1}{\sqrt{2\pi b}}
\left[s-\frac{2s^3-s}{12b}\right]
\left[1+O\left(\frac{1}{b^2\sqrt{b}}\right)\right]\right).
$$
\label{expect}
\end{lemma}

{\it Proof.}
Compute
$$
{\sf E}\left[\eta_b^{(s)}I(\eta_b<b)\right]=\sum_{i=0}^{b-1}\frac{i!}{(i-s)!}\frac{b^i}{i!}e^{-b}=
$$
$$
=b^s\sum_{i=s}^{b-1}\frac{b^{i-s}}{(i-s)!}e^{-b}=b^s{\sf P}(\eta_b<b-s)=b^s\left[{\sf P}(\eta_b<b)-{\sf P}(b-s\leq\eta_b<b)\right].
$$
By the Stirling's approximation,
$$
{\sf P}(b-s\leq\eta_b<b)=\sum_{i=1}^{s}\frac{b^{b-i} e^{-b}}{(b-i)!}=
\sum_{i=1}^{s}\frac{1}{\sqrt{2\pi b}} e^{-\frac{i^2}{2(b-i)}-\frac{1}{12(b-i)}}\left(1+\frac{i}{2(b-i)}\right)\left(1+O\left(\frac{1}{b^2}\right)\right)=
$$
$$
\frac{1}{\sqrt{2\pi b}}\left[\sum_{i=1}^{s}\left(1-\frac{i^2}{2b}-\frac{1}{12b}+\frac{i}{2b}\right)\right]
\left(1+O\left(\frac{1}{b^2}\right)\right)
=\frac{1}{\sqrt{2\pi b}}\left[s-\frac{2s^3-s}{12b}\right]
\left(1+O\left(\frac{1}{b^2}\right)\right).\quad\Box
$$

Since, for $s,k\in\mathbb{N}$, $s<k$, 
$
\sum_{i=0}^k (-1)^i {k\choose i} i^{(s)}=0,
$
Lemma~\ref{expect} with (\ref{p_b_expantion}) imply the following.

\begin{claim}
For every $k\in\mathbb{N}$, set $h_k^1={\sf E}\left[(b-\eta_b)^k I(\eta_b<b)\right], h_k^2={\sf E}\left[\eta_b(b-\eta_b)^k I(\eta_b<b)\right].$
Then 
$$
h_1^1=\frac{\sqrt{b}}{\sqrt{2\pi}}-\frac{1}{12\sqrt{2\pi b}}+O(b^{-1}),\quad h_2^1=b\left(\frac{1}{2}-\frac{1}{3\sqrt{2\pi b}}\right)+O(1),\quad h_3^1=\frac{2b\sqrt{b}}{\sqrt{2\pi}}+O(b),
$$
$$
h_k^1=O(b^{k-2})\text{ for all }k\geq 4;\quad h_1^2=\frac{b\sqrt{b}}{\sqrt{2\pi}}+O(b),\quad h_k^2=O(b^{k})\text{ for all }k\geq 2.
$$
\label{expectations2}
\end{claim}

Using Claim~\ref{expectations2}, we get the following asymptotical expansion of $z_{b,n}$ when $b=O(\sqrt{n})$.

\begin{claim}
Let $C>0$ and $n>Cb^2$. Then $z_{b,n}=\frac{1}{3}+\frac{4}{135b}+\frac{b}{3n}+O\left(\frac{1}{b\sqrt{b}}\right).$
\label{CL_z_estimation}
\end{claim}
{\it Proof.} With Stirling’s approximation, we have 
$$
 p_{b,n}= \left(1+O(b^{-2})\right)\left({\sf P}(\eta_b<b)+\frac{1}{2n}{\sf E}\eta_b I(\eta_b<b)-\frac{1}{2n}h_2^1-\frac{1}{3n^2}h_3^1\right).
$$
Applying Lemma 1 and Claim 4 and substituting the final asymptotics into the definition~(\ref{z_bn_def}) of $z_{b,n}$, we get the result. 
See the details in Appendix A.

\section{Proof of Theorem~\ref{ps}}
\label{th2}

For boundary cases, when $b\leq 5$ or $b\geq n - 5$, the proof is given in Appendix B. Below, we prove Theorem~\ref{ps} for $6\leq b\leq n-6$.\\

Denote 
\begin{equation}
g_{b,n}:=\int_{1-(b+1)/n}^{1-b/n}(1-z)^{b-1}z^{n-b}dz.
\label{g_bn_defintion}
\end{equation}
From  Claim~\ref{CL_below_0.5}, we get
\begin{equation}
\begin{aligned}
 \sum_{\ell=0}^3\frac{1}{(\ell+1)!n^{\ell+1}}\frac{\partial^{\ell} g}{\partial z^{\ell}}\left(1-\frac{b+1}{n}\right)+
 \frac{1}{120n^{5}}\frac{\partial^{4} g}{\partial z^{4}}\left(1-\frac{b+1}{n}\right)=:g^-_{b,n}\leq\\
 g_{b,n}\leq g^+_{b,n}:=\sum_{\ell=0}^3\frac{1}{(\ell+1)!n^{\ell+1}}\frac{\partial^{\ell} g}{\partial z^{\ell}}\left(1-\frac{b+1}{n}\right)+
 \frac{1}{120n^{5}}\frac{\partial^{4} g}{\partial z^{4}}\left(1-\frac{b}{n}\right).
\end{aligned}
\label{g_bounds}
\end{equation}

Denote 
\begin{equation}
\begin{split}
P_{b,n}:=&12b^5-16b^4n+64b^4+4b^3n^2-71b^3n+138b^3+16b^2n^2-112b^2n+ \\ 
&156b^2+12bn^2-105bn+94b+24n^2-48n+24.
\end{split}
\label{Pbn_def}
\end{equation}

Clearly,
$$
\left(\frac{b+1}{n}\right)^b\left(1-\frac{b+1}{n}\right)^{n-b}-bg_{b,n}^{\gamma} =
$$
\begin{equation}
\frac{1}{120n^6}\left(5\left(\frac{b+1}{n}\right)^{b-4}\left(1-\frac{b+1}{n}\right)^{n-b-2}P_{b,n}-bn\frac{\partial^{4} g}{\partial z^{4}}\left(1-\frac{b+r(\gamma)}{n}\right)\right).\\
\label{g_equality}
\end{equation}

where $\gamma\in\{+,-\}$; $r=0$ when $\gamma=+$ and $r=1$ when $\gamma=-$.\\

Below, we will use the following bounds proven in Appendix C.

\begin{claim}
If $6\leq b\leq\frac{n-2}{3}$, then
\begin{equation}
\label{small_b_ineq}
\frac{\partial^{4} g}{\partial z^{4}}\left(1-\frac{b}{n}\right)<\frac{5}{bn}\left(\frac{b+1}{n}\right)^{b-4}\left(1-\frac{b+1}{n}\right)^{n-b-2}P_{b,n}.
\end{equation}
If $5\leq b<\frac{n}{3}$, then
$$
\frac{\partial^{4} g}{\partial z^{4}}\left(1-\frac{b+1}{n}\right)>\frac{5}{bn}
\left(\frac{b+1}{n}\right)^{b-4}\left(1-\frac{b+1}{n}\right)^{n-b-2}P_{b,n}
$$
\begin{equation}
+\frac{120n^5}{b}\left(\frac{n-b}{n-b-1}-\left(\frac{b}{b+1}\right)^b\left(\frac{n-b}{n-b-1}\right)^{n-b}\right)\left(\frac{b+1}{n}\right)^b\left(1-\frac{b+1}{n}\right)^{n-b}.
\label{negativeness2}
\end{equation}

If $\frac{n-1}{3}\leq b\leq\frac{n}{2}$, then
\begin{equation}
\frac{\partial^{4} g}{\partial z^{4}}\left(1-\frac{b+1}{n}\right)>\frac{5}{bn}
\left(\frac{b+1}{n}\right)^{b-4}\left(1-\frac{b+1}{n}\right)^{n-b-2}P_{b,n}.
\label{negativeness}
\end{equation}
\label{th1_main_claim}
\end{claim}


Further, we distinguish three cases: $6\leq b\leq (n-2)/3$,  $(n-1)/3\leq b \leq n/2$ and $n - 5 > b>n/2$.\\

1) Assume that $6\leq b\leq (n-2)/3$.


Due to~(\ref{p_difference_g})~and~(\ref{g_bounds}), it is sufficient to prove that $\left(\frac{b+1}{n}\right)^b\left(1-\frac{b+1}{n}\right)^{n-b}-bg^+_{b,n}>0$. The latter inequality follows from~(\ref{g_equality}) and Claim~\ref{th1_main_claim}.\\

2) Let $(n-1)/3\leq b \leq  n/2$.

Due to~(\ref{p_difference_g})~and~(\ref{g_bounds}), it is sufficient to prove that $\left(\frac{b+1}{n}\right)^b\left(1-\frac{b+1}{n}\right)^{n-b}-bg^-_{b,n}<0$. The latter inequality follows from~(\ref{g_equality}) and Claim~\ref{th1_main_claim}.\\

3) Finally, let $n/2<b<n-5$. By the definition,
\begin{equation}
 p_{b,n}=1-{n\choose b}\left(\frac{b}{n}\right)^b\left(1-\frac{b}{n}\right)^{n-b}-p_{n-b,n}.
\label{before_above}
\end{equation}
Thus,
\begin{equation}
p_{b+1,n}-p_{b,n}=\frac{(n-1)!}{b!(n-b-1)!}\frac{(b+1)^b(n-b-1)^{n-b-1}}{n^{n-1}}B+(p_{n-b,n}-p_{n-b-1,n}).
\label{p_above_n_over_2}
\end{equation}
where $B=\left(1-\frac{1}{b+1}\right)^b\left(1+\frac{1}{n-b-1}\right)^{n-b-1}-1$. For all $b\geq (n+1)/2$ and $n\geq6$, we have 
$B<\exp\left(\frac{1}{2(b+1)}+\frac{1}{6(b+1)^2}-\frac{1}{2(n-b-1)}+\frac{1}{3(n-b-1)^2}\right)-1<0
$, which proves $p_{b+1,n} < p_{b,n}$ for $(n+1)/2\leq b< 2n/3$.
\\
Let $n-6\geq b\geq 2n/3$. Set $\tilde b=n-b-1$. In this case $5\leq \tilde b < n/3$. From Claim~\ref{p_difference_g}, Claim~\ref{th1_main_claim} and relations~(\ref{g_bounds}),~(\ref{g_equality})~and~(\ref{p_above_n_over_2}), we get that $p_{b+1,n}-p_{b,n}=$
$$
{n\choose\tilde b}
\left[\left(1+\frac{n-\tilde b}{n-\tilde b-1}-\left(\frac{\tilde b}{\tilde b+1}\right)^{\tilde b}\left(\frac{n-\tilde b}{n-\tilde b-1}\right)^{n-\tilde b}\right)\left(\frac{\tilde b+1}{n}\right)^{\tilde b}\left(1-\frac{\tilde b+1}{n}\right)^{n-\tilde b}
-\tilde bg_{\tilde b,n}\right]<
$$
$$
\left(\frac{\tilde b + 1}{n}\right)^{\tilde b-4}\left(1-\frac{\tilde b + 1}{n}\right)^{n-\tilde b - 2}{n\choose\tilde b}\times\biggl[\frac{1}{24n^6}P_{\tilde b,n}-\frac{\tilde b\frac{\partial^{4} g}{\partial z^{4}}\left(1-\frac{\tilde b+1}{n}\right)}{120n^5\left(\frac{\tilde b+1}{n}\right)^{\tilde b}\left(1-\frac{\tilde b+1}{n}\right)^{n-\tilde b-2}} + 
$$
$$
\left(\frac{n-\tilde b}{n-\tilde b-1} - \left(\frac{\tilde b}{\tilde b+1}\right)^{\tilde b}\left(\frac{n-\tilde b}{n-\tilde b-1}\right)^{n-\tilde b}\right)\left(\frac{\tilde b + 1}{n}\right)^4\left(1-\frac{\tilde b + 1}{n}\right)^2\biggr]<0.
$$


\section{Proof of Theorem~\ref{zs}}
\label{th1}

If $b=(n-1)/2$, then, from Theorem~\ref{ps}, for $n$ large enough,
$$
 \frac{z_{b+1,n}}{z_{b,n}}=\frac{1/2-p_{b+1,n}}{1/2-p_{b,n}}\frac{b^b(n-b)^{n-b}}{(b+1)^{b+1}(n-b-1)^{n-b-1}}=\frac{1/2-p_{b+1,n}}{1/2-p_{b,n}}>1.
$$

Let $b<\frac{n-1}{2}$ (we consider the case $b\geq n/2$ in the very end of the proof). In this case,  \begin{equation}
b\ln\left(1+\frac{1}{b}\right)+(n-b-1)\ln\left(1-\frac{1}{n-b}\right)<-\frac{n-2b}{2b(n-b)}+\frac{1}{2(n-b)^2}+\frac{1}{3b^2}<0.
\label{negative_appendix}
\end{equation}

Since
$$
\left[\int_{1-(b+1)/n}^1-\int_0^{1-(b+1)/n}\right]g(z)dz-\left[\int_{1-b/n}^1-\int_0^{1-b/n}\right]g(z)dz=2\int_{1-(b+1)/n}^{1-b/n}g(z)dz,
$$
by~(\ref{g_useful})~and~(\ref{g_bn_defintion}), 
$$
z_{b+1,n}-z_{b,n}\geq\frac{n^{n}}{(n-b)(b+1)^b(n-b-1)^{n-b-1}} \times \biggl[
 bg_{b,n}- \left(\frac{b+1}{n}\right)^b\times
 $$
 $$
 \left(\frac{n-b-1}{n}\right)^{n-b}+
 \frac{b}{2}\left(1-\left(1+\frac{1}{b}\right)^b\left(1-\frac{1}{n-b}\right)^{n-b-1}\right)\left[\int_{1-b/n}^1-\int_0^{1-b/n}\right]g(z)dz\biggr].
$$
Denote 
$$
A_{b,n}:=(n-b)\left(\frac{b+1}{n}\right)^5\left(1-\frac{b+1}{n}\right)^3,\quad x=\frac{b+1}{n}.
$$
Using~(\ref{g_bounds})~and~(\ref{g_equality}), we get
\begin{equation}
\begin{aligned}    
 A_{b,n}(z_{b+1,n}-z_{b,n})>\frac{b\frac{\partial^{4} g}{\partial z^{4}}\left(1-x\right)}{120n^4 x^{b-5}\left(1-x\right)^{n-b-4}}-\frac{x(1-x)^2}{24n^5}P_{b,n}+\left(\frac{n}{b}\right)^b\left(\frac{n}{n-b}\right)^{n-b}\times\\
 x^5\left(1-x\right)^4\frac{nb}{2}\left(\left(\frac{b}{b+1}\right)^b\left(\frac{n-b}{n-b-1}\right)^{n-b}-\frac{n-b}{n-b-1}\right)\left[\int_{1-b/n}^1-\int_0^{1-b/n}\right]g(z)dz.
\end{aligned}
\label{delta_z_lower}
\end{equation}
Since, for all $y>0$, $e^y>1+y+y^2/2$ and, for all $y>0$ small enough, $\ln(1-y)>-y-y^2/2-y^3/3-y^4/3$ and $\ln(1+y)>y-y^2/2$, we get that, for all $b$ and $n-b$ large enough,
$$
 \left(\frac{b}{b+1}\right)^b\left(\frac{n-b}{n-b-1}\right)^{n-b}-\frac{n-b}{n-b-1}>
$$
$$
 \exp\left[b\left(-\frac{1}{b+1}-\frac{1}{2(b+1)^2}-\frac{1}{3(b+1)^3}-\frac{1}{3(b+1)^4}\right)+(n-b)\left(\frac{1}{n-b-1}-\frac{1}{2(n-b-1)^2}\right)\right]
$$
$$
 -\frac{n-b}{n-b-1}>\exp\left[\frac{1}{2(b+1)}+\frac{1}{6(b+1)^2}+\frac{1}{2(n-b-1)}-\frac{1}{2(n-b-1)^2}\right]-1-\frac{1}{n-b-1}>
$$
\begin{equation}
 \frac{1}{2(b+1)}+\frac{1}{6(b+1)^2}-\frac{1}{2(n-b-1)}-\frac{1}{2(n-b-1)^2}+\frac{1}{8(b+1)^2}.
\label{exp_lower}
\end{equation}

The rest is divided into 3 parts: 1) $(n-1)/2\geq b>(1+\varepsilon)\sqrt{\frac{31}{60}n}$; 2)  $(1+\varepsilon)\sqrt{\frac{31}{60}n}\geq b>(1+\varepsilon)\sqrt{\frac{4}{45}n}$; 3) $b<(1-\varepsilon)\sqrt{\frac{4}{45}n}$.\\

1) Let $(n-1)/2>b>(1+\varepsilon)\sqrt{\frac{31}{60}n}$. From Theorem~\ref{Sam} and (\ref{zs_integral}),
\begin{equation}
\label{gz_integral_lowerbound}
 \left[\int_{1-b/n}^1-\int_{0}^{1-b/n}\right]g(z)dz\geq\frac{2}{3b}(b/n)^b(1-b/n)^{n-b}.
\end{equation}
By~(\ref{delta_z_lower})~and~(\ref{exp_lower}), we get that, for $n$ large enough,
$$
 A_{b,n}(z_{b+1,n}-z_{b,n})>\frac{b\frac{\partial^{4} g}{\partial z^{4}}\left(1-x\right)}{120n^4 x^{b-5}\left(1-x\right)^{n-b-4}}-\frac{x(1-x)^2}{24n^5}P_{b,n}+
$$
$$ 
 x^4\left(1-x\right)^{4} \left(\frac{1}{6}+\frac{1}{18(b+1)}\right)-
 x^5(1-x)^{3}\left(\frac{1}{6}+\frac{1}{6(n-b-1)}\right),
$$
It is easy to see that, for every $\beta>0$, there exists $\delta>0$ such that, for all $n$ large enough and $b>(1+\varepsilon)\sqrt{\frac{31}{60}n}$, 
$$
 A_{b,n}(z_{b+1,n}-z_{b,n})>(1-\beta)\frac{1}{6}x^5(1-x)^3\geq 0,\quad\text{if }x\geq\delta,
$$
$$
 A_{b,n}(z_{b+1,n}-z_{b,n})>(1-\beta)\left[\frac{1}{6}x^5-\frac{31}{360}\frac{x^3}{n}\right]\geq 0,\quad\text{if }x<\delta.
$$

2) Let $(1+\varepsilon)\sqrt{\frac{31}{60}n}\geq b>(1+\varepsilon)\sqrt{\frac{4}{45}n}$. By Claim~\ref{CL_z_estimation}, for every $\gamma>0$ and $n$ large enough, $z_{b,n}\geq\frac{1}{3}+\left[\frac{4}{135b}+\frac{b}{3n}\right](1-\gamma)$.
Thus, for such $n$, by (\ref{zs_integral}),
$$
 \left[\int_{1-b/n}^1-\int_{0}^{1-b/n}\right]g(z)dz\geq\left(\frac{2}{3b}+\left[\frac{8}{135b^2}+\frac{2}{3n}\right](1-\gamma)\right)(b/n)^b(1-b/n)^{n-b}.
$$
By~(\ref{delta_z_lower})~and~(\ref{exp_lower}), we get that,  for every $\tilde\gamma>\gamma$ and $n$ large enough,
$$
 A_{b,n}(z_{b+1,n}-z_{b,n})>\frac{b\frac{\partial^{4} g}{\partial z^{4}}\left(1-x\right)}{120n^4 x^{b-5}\left(1-x\right)^{n-b-4}}-\frac{x(1-x)^2}{24n^5}P_{b,n}+
$$
$$ 
 x^4\left(1-x\right)^{4} \left(\frac{1}{6}+\left(\frac{7}{72}+\frac{2}{135}(1-\tilde\gamma)\right)\frac{1}{b+1}\right)+
 \frac{1}{6}x^5\left(1-x\right)^{4}(1-\tilde\gamma)
 -\frac{1}{6}x^5\left(1-x\right)^{3}.
$$
Therefore, for every $\beta>0$, large enough $n$ and $(1+\varepsilon)\sqrt{\frac{31}{60}n}\geq b>(1+\varepsilon)\sqrt{\frac{4}{45}n}$,
\begin{equation}
A_{b,n}(z_{b+1,n}-z_{b,n})>(1-\beta)\left[\frac{x^5}{3}-\frac{4x^3}{135n}\right]>0.
\label{convenient_inequality}
\end{equation}

3) Let $b<(1-\varepsilon)\sqrt{\frac{4}{45}n}$. For constant $b$ and large enough $n$, $z_{b+1,n}<z_{b,n}$ since $y_{b+1}<y_b$. Below, we consider $b$ as large as desired. By Claim~\ref{CL_z_estimation}~and~(\ref{zs_integral}), for every $\gamma>0$ and $n$ large enough,
$$
 \left[\int_{1-b/n}^1-\int_{0}^{1-b/n}\right]g(z)dz\leq\left(\frac{2}{3b}+\left[\frac{8}{135b^2}+\frac{2}{3n}\right](1+\gamma)\right)(b/n)^b(1-b/n)^{n-b}.
$$
Note that, here, we should use $g^+_{b,n}$ instead of $g^-_{b,n}$, but the contribution of the difference between them in (\ref{convenient_inequality}) is at most $O(x^2/n^2)$, and so, for some constant $c$ and every positive $\beta$,
$$
A_{b,n}(z_{b+1,n}-z_{b,n})<(1+\beta)\left[\frac{x^5}{3}-\frac{4x^3}{135n}+c\frac{x^2}{n^2}\right]<0
$$
for $n$ large enough.\\


Finally, let us consider $b\geq n/2$. From~(\ref{before_above}), we get
$$
z_{b,n}=1-\frac{1/2-\left(1-{n\choose b}(b/n)^b(1-b/n)^{n-b}-p_{n-b,n}\right)}{{n\choose b}(b/n)^b(1-b/n)^{n-b}}=1-z_{n-b,n}.
$$
Therefore, $z_{b+1,n}>z_{b,n}$ if and only if $z_{n-b-1,n}<z_{n-b,n}$. Theorem is proved.

\section{Discussions}
\label{further}

It is natural to ask, do our results remain true for binomial distributions having other second parameters but still approaching Poisson distribution? Given $c>0$, a possible generalization of Theorem~\ref{ps} is a study of behaviour of $p_{c,b,n}={\sf P}(\xi_{c,b,n}<b)$ as a function of $b$, where $\xi_{c,b,n}$ has binomial distribution with parameters $n$ and $\frac{b}{n+c}$. Such a study is highly related to the open {\it problem of small deviations inequality}~(see, e.g.,~\cite{Feige}), which can be formulated as: for $c>0$, find the minimum of ${\sf P} \left(X_1 + X_2 + \ldots + X_n < n + c\right)$ over all sets of independent non-negative unit mean random variables $\{X_i\}_{i=1}^n$. It was shown~(see, e.g.~\cite{He}) that the optimal random variables are two-point. If we further restrict our set to independent identically distributed random variables, we will reduce (see the details below) the initial problem to the analysis of monotonicity of $p_{c,b,n}$. More formally, for every $\mathcal{X}\in[0,1)\times(1,\infty)$, consider the distribution $\mathcal{Q}(\mathcal{X})$ on $\mathcal{X}$ with mean 1. Let  $X_1,\ldots,X_n$ be independent random variables with identical distribution $\mathcal{Q}(X)$. The problem is to find $\min_{\mathcal{X}\in[0,1)\times(1,\infty)}{\sf P}(X_1+\ldots+X_n<n+c)$. In what follows, for simplicity of notations and computations, we consider $c=1$ (for arbitrary $c>0$, similar arguments work). In~\cite{Samuels_sum}, S. M. Samuels proved that the minimum equals $\left(\frac{n}{n+1}\right)^n$. Let us show that this is true if and only if, for all integers $b\geq 2$, $n\geq 2b$, $p_{1,1,n}\leq p_{1,b,n}$.

If $\min_{\mathcal{X}\in[0,1)\times(1,\infty)}{\sf P}(X_1+\ldots+X_n<n+1)=\left(\frac{n}{n+1}\right)^n$, then for $b\geq 2$, $n\geq 2b$, consider $\mathcal{X}=\{0,\frac{n+1}{b}\}$. Clearly, ${\sf P}(X_1+\ldots+X_n<n+1)=p_{1,b,n}\geq p_{1,1,n}=\left(\frac{n}{n+1}\right)^n$. Now, assume that, for all integers $b\geq 2$, $n\geq 2b$, $p_{1,1,n}\leq p_{1,b,n}$. Let $(\alpha,\beta)=\mathcal{X}\in[0,1)\times(1,\infty)$. Then ${\sf P}(X_1=\beta)=\frac{1-\alpha}{\beta-\alpha}$. 
Then, clearly, ${\sf P}(X_1+\ldots+X_n<n+1)={\sf P}(\zeta<\lceil\frac{n+1-n\alpha}{\beta-\alpha}\rceil)$, where $\zeta\sim\mathrm{Bin}(n,\frac{1-\alpha}{\beta-\alpha})$. Denote $b:=\lceil\frac{n+1-n\alpha}{\beta-\alpha}\rceil$. It remains to prove that ${\sf P}(\zeta<b)\geq p_{1,b,n}$. The latter inequality follows from the fact that $\frac{1-\alpha}{\beta-\alpha}\leq\frac{b}{n+1}$. Indeed, $\frac{b}{n+1}=
\left\lceil\frac{n+1-n\alpha}{\beta-\alpha}\right\rceil\frac{1}{n+1}\geq \frac{1-\frac{n}{n+1}\alpha}{\beta-\alpha}\geq\frac{1-\alpha}{\beta-\alpha}$.

The inequality $p_{1,1,n}\leq p_{1,b,n}$ can be derived from a result on a monotonicity of $p_{1,b,n}$ analogous to Theorem~\ref{ps}. Notice that, by following the proof of Claim~\ref{pg_expression}, we get that the monotonicity of $p_{1,b,n}$ boils down to the analysis of the same function $g(z)$ introduced in the beginning of Section~\ref{tools}, but now this function must be studied on $[1-(b+1)/(n+1),1-b/(n+1)]$. We conjecture that $p_{1,b,n}$ increases in $b$  (see Appendix D for empirical comparison of $p_{b,n}$ and $p_{c,b,n}$). Notice that this conjecture imply even more than the result of S.M. Samuels. 
If $p_{1,b,n}$ increases, we get
\begin{conjecture}
Let $X=(\alpha,\beta)\in[0,1)\times(1,\infty)$. Let $b$ be an integer and $ b-1<\frac{n+1-n\alpha}{\beta-\alpha}\leq b.$
Then ${\sf P}(\xi_1+\ldots+\xi_n<n+1)\geq p_{1, b,n}$, where $\xi_1,\ldots,\xi_n$ are independent random variables with identical distribution $\mathcal{Q}(X)$, and the equality holds if and only if $\alpha=0$ and $\frac{n+1}{\beta}=b$.
\end{conjecture}


\section*{Acknowledgements} The work is supported by Mathematical Center in Akademgorodok, the agreement with Ministry of Science and High Education of the Russian Federation number  075-15-2019-1613.

\section*{Appendix A. Missing proofs of Section~\ref{tools}}

\subsection*{Proof of Claim~\ref{pg_expression}}

Let us rewrite the definition of $p_{b,n}$ in the following way:
$$
p_{b,n}=\sum_{i=0}^{b-1}{n\choose b-1}{b-1\choose i}\frac{(n-b+1)!(b-1-i)!}{(n-i)!}\left(\frac{b}{n}\right)^i\left(1-\frac{b}{n}\right)^{n-i}.
$$
Since 
$$
\frac{(n-b)!(b-1-i)!}{(n-i)!}=\frac{\Gamma(n-b+1)\Gamma(b-i)}{\Gamma(n-i+1)}=\B(n-b+1,b-i)=\int_0^1 x^{n-b}(1-x)^{b-i-1}dx,
$$
we get
$$
p_{b,n}=
n{n-1\choose b-1}\times
$$
$$
\int_0^1\left[\sum_{i=0}^{b-1}{b-1\choose i}(1-x)^{b-i-1}\left(\frac{b}{n}\right)^i\left(1-\frac{b}{n}\right)^{b-i-1}\right]x^{n-b}\left(1-\frac{b}{n}\right)^{n-b+1}dx=
$$
$$
\frac{(n-1)!}{(n-b-1)!(b-1)!}\int_0^1\left[\frac{b}{n}+(1-x)\left(1-\frac{b}{n}\right)\right]^{b-1}x^{n-b}\left(1-\frac{b}{n}\right)^{n-b}dx=
$$
$$
\frac{n!}{(n-b)!(b-1)!}\int_0^1\left[1-x\left(1-\frac{b}{n}\right)\right]^{b-1}\left[x\left(1-\frac{b}{n}\right)\right]^{n-b}d\left[x\left(1-\frac{b}{n}\right)\right]=
$$
$$
\frac{\Gamma(n+1)}{\Gamma(n-b+1)\Gamma(b)}\int_0^{1-b/n}(1-z)^{b-1}z^{n-b}dz=
\frac{\int_0^{1-b/n}(1-z)^{b-1}z^{n-b}dz}{\int_0^1(1-z)^{b-1}z^{n-b}dz}.
$$
Therefore,
$$
 p_{b+1,n}-p_{b,n}=
 \frac{\int_0^{1-(b+1)/n}(1-z)^b z^{n-b-1}dz}{\int_0^1(1-z)^{b}z^{n-b-1}dz}-\frac{\int_0^{1-b/n}(1-z)^{b-1}z^{n-b}dz}{\int_0^1(1-z)^{b-1}z^{n-b}dz}=
$$
$$
 \frac{\int_0^{1-(b+1)/n}(1-z)^b d\left(z^{n-b}\right)}{\int_0^1(1-z)^{b}d\left(z^{n-b}\right)}-\frac{\int_0^{1-b/n}(1-z)^{b-1}z^{n-b}dz}{\int_0^1(1-z)^{b-1}z^{n-b}dz}=
$$
$$
 \frac{\left(\frac{b+1}{n}\right)^b\left(1-\frac{b+1}{n}\right)^{n-b}+b\int_0^{1-(b+1)/n}(1-z)^{b-1} z^{n-b}dz}{b\int_0^1(1-z)^{b-1}z^{n-b}dz}-\frac{\int_0^{1-b/n}(1-z)^{b-1}z^{n-b}dz}{\int_0^1(1-z)^{b-1}z^{n-b}dz}=
$$
$$ 
\frac{\left(\frac{b+1}{n}\right)^b\left(1-\frac{b+1}{n}\right)^{n-b}-b\int_{1-(b+1)/n}^{1-b/n}(1-z)^{b-1} z^{n-b}dz}{b\int_0^1(1-z)^{b-1}z^{n-b}dz}.
$$

Now let us analyze the values $z_{b,n}$. Using~(\ref{from_p_to_beta}), we get
$$
z_{b,n} =\frac{\frac{1}{2}b\left(\left[\int_{1-b/n}^1-\int_{0}^{1-b/n}\right]g(z)dz\right)}{(b/n)^b(1-b/n)^{n-b}}.
$$
Then
$$
z_{b+1,n}-z_{b,n}=
 \frac{n^{n}}{2}\biggl(\frac{\left[\int_{1-(b+1)/n}^1-\int_{0}^{1-(b+1)/n}\right](1-z)^bd(z^{n-b})}{(b+1)^b(n-b-1)^{n-b-1}(n-b)}-
$$
$$ 
 \frac{\left[\int_{1-b/n}^1-\int_{0}^{1-b/n}\right]g(z)dz}{b^{b-1}(n-b)^{n-b}}\biggr)=\frac{n^{n}}{2(n-b)}\times
 $$
 $$
 \biggl(\frac{-2\left(\frac{b+1}{n}\right)^b\left(\frac{n-b-1}{n}\right)^{n-b}+ 
 b\left[\int_{1-(b+1)/n}^1-\int_{0}^{1-(b+1)/n}\right]g(z)dz}{(b+1)^b(n-b-1)^{n-b-1}}- \frac{b\left[\int_{1-b/n}^1-\int_{0}^{1-b/n}\right]g(z)dz}{b^b(n-b)^{n-b-1}}\biggr)
$$
$$
=\frac{n^{n}}{2(n-b)(b+1)^b(n-b-1)^{n-b-1}}\times
\biggl(-2\left(\frac{b+1}{n}\right)^b\left(\frac{n-b-1}{n}\right)^{n-b}+
$$
$$
b\left[\int_{1-(b+1)/n}^1-\int_{0}^{1-(b+1)/n}\right]g(z)dz-
$$
$$
b\left(1+\frac{1}{b}\right)^b\left(1-\frac{1}{n-b}\right)^{n-b-1}\left[\int_{1-b/n}^1-\int_{0}^{1-b/n}\right]g(z)dz\biggr).
$$
$\quad\Box$

\subsection*{Proof of Claim~\ref{CL_below_0.5}}
For every $\ell\in\{1,\ldots,\min\{b-1,n-b\}\}$, denote 
$$
f_{\ell}=\frac{\partial^{\ell} g/\partial z^{\ell}}{(1-z)^{b-1-\ell}z^{n-b-\ell}}.
$$
It is easy to see that, for every $\ell$, $\partial f_{\ell+1}/\partial z=-(\ell+1)(n-\ell-1)f_{\ell}$, and
$$
 f_2(z)=z^2(n-1)(n-2)-2z(n-2)(n-b)+(n-b)(n-b-1)
$$
is negative on 
$$
\Upsilon:=\left(\frac{n-b}{n-1}-\frac{1}{n-1}\sqrt{\frac{(n-b)(b-1)}{n-2}},\frac{n-b}{n-1}+\frac{1}{n-1}\sqrt{\frac{(n-b)(b-1)}{n-2}}\right).
$$
Let us show that $\Delta_{b,n}\subset\Upsilon$. First, 
$$
1-\frac{b+1}{n}>\frac{n-b}{n-1}-\frac{1}{n-1}\sqrt{\frac{(n-b)(b-1)}{n-2}}
$$ 
since the difference between the left side and right side of this inequality equals $\sqrt{\frac{(n-b)(b-1)}{n-2}}-\frac{2n-b-1}{n}$, and
$$
 n^2(n-b)(b-1)-(2n-b-1)^2(n-2)=
$$
$$
n(b-5)(n(n-b)-b)+n(12n-15b-9)+2b^2+2+4b>0
$$
because $5\leq b\leq n/2$. Second, 
$$
1-\frac{b}{n}<\frac{n-b}{n-1}+\frac{1}{n-1}\sqrt{\frac{(n-b)(b-1)}{n-2}}
$$ 
since, in fact, $1-\frac{b}{n}<\frac{n-b}{n-1}$. So, $f_3(z)$ increases on $\Delta_{b,n}$. 

Now, let us show that 
$$
f_3(1-\frac{b}{n})=(n-b)\left[\left(2-7\frac{n-b}{n}\right)\left(1-\frac{n-b}{n}\right)+6\frac{(n-b)^2}{n^3}\right]<0.
$$
The derivative of the function $(2-7(x/n))(1-(x/n))+6(x^2/n^3)$ with respect to $x$ is negative when $x<\frac{9n^2}{14n+12}$ and positive when $x>\frac{9n^2}{14n+12}$. Therefore, negativeness of this function in $x=\frac{2}{7}n+1$ and $x=n-2$ implies its negativeness for all $x\in\left[\frac{2}{7}n+1,n-2\right]$.

As $f_3$ increases on $\Delta_{b,n}$ and is negative in $1-b/n$, $f_4^{\prime}(z)$ is positive on $\Delta_{b,n}$. Moreover, the derivative of $(1-z)^{b-5}z^{n-b-4}$ with respect to $z$ is also positive on $\Delta_{b,n}$. Thus, $\partial^4 g/\partial z^4$ increases on $\Delta_{b,n}$, and this implies the statement of Claim~\ref{CL_below_0.5}. $\quad\Box$

\subsection*{Proof of Claim~\ref{CL_z_estimation}}
First, let us estimate $p_b-p_{b,n}$. Let $n>Cb^2$. Then
$$
 p_{b,n}=\sum_{i=0}^{b-1}{n\choose i}\left(\frac{b}{n}\right)^i\left(1-\frac{b}{n}\right)^{n-i}=
 $$
 $$
 \left(1+O(b^{-2})\right)\sum_{i=0}^{b-1}\frac{1}{i!e^i(n-i)^{n-i}}\sqrt{\frac{n}{n-i}} b^i(n-b)^{n-i}=
$$
$$ 
 \left(1+O(b^{-2})\right)\sum_{i=0}^{b-1}\frac{b^i}{i!e^b}\left(1+\frac{i}{2n}\right) e^{-\frac{(b-i)^2}{2(n-i)}-\frac{(b-i)^3}{3(n-i)^2}}=
$$ 
$$ 
 \left(1+O(b^{-2})\right)\sum_{i=0}^{b-1}\frac{b^i}{i!e^b}\left(1+\frac{i}{2n}\right)
 e^{-\frac{(b-i)^2}{2n}\left(1+\frac{i}{n}\right)}\left(1-\frac{(b-i)^3}{3 n^2}\right).
$$
Since, for any positive $x$ and $y$, 
$$
1-x(1+y)<e^{-x(1+y)}<e^{-x}<1-x+\frac{x^2}{2},
$$
by Claim~\ref{expectations2}, we get
$$
 p_{b,n}=\left(1+O(b^{-2})\right)\sum_{i=0}^{b-1}\frac{b^i}{i!e^b}\left(1+\frac{i}{2n}-\frac{(b-i)^2}{2n}-\frac{(b-i)^3}{3 n^2}\right)=
$$
$$ 
 \left(1+O(b^{-2})\right)\left({\sf P}(\eta_b<b)+\frac{1}{2n}{\sf E}\eta_b I(\eta_b<b)-\frac{1}{2n}h_2^1-\frac{1}{3n^2}h_3^1\right).
$$
By (\ref{p_b_expantion}), Lemma~\ref{expect} and Claim~\ref{expectations2}, we get
$$ 
 p_{b,n}=\frac{1}{2}-\frac{1}{3\sqrt{2\pi b}}-\frac{1}{540\sqrt{2\pi}b\sqrt{b}}-\frac{\sqrt{b}}{2n\sqrt{2\pi}}+O(b^{-2}).
$$
Therefore,
$$
 z_{b,n}=\frac{\frac{1}{2}-p_{b,n}}{{n\choose b}\left(\frac{b}{n}\right)^b\left(1-\frac{b}{n}\right)^{n-b}}=
 \sqrt{\frac{n-b}{n}}e^{\frac{1}{12}b}\left(\frac{1}{3}+\frac{1}{540b}+\frac{b}{2n}+O\left(\frac{1}{b\sqrt{b}}\right)\right)=
$$
$$
\left(1-\frac{b}{2n}\right)\left(1+\frac{1}{12}b\right)\left(\frac{1}{3}+\frac{1}{540b}+\frac{b}{2n}\right)+O\left(\frac{1}{b\sqrt{b}}\right)=
\frac{1}{3}+\frac{4}{135b}+\frac{b}{3n}+O\left(\frac{1}{b\sqrt{b}}\right).\quad\Box
$$

\section*{Appendix B. Boundary cases of Theorem~\ref{ps} ($b \leq 5$ and $b \geq n - 5$)}

\begin{claim}
Let $b \leq 5$. Then

\begin{center}

$p_{b,n} < p_{b+1, n}$ if $n \geq 3b+2$; 

$p_{b,n} > p_{b+1, n}$ if $b + 1 \leq n\leq 3b+1$.

\end{center}

\label{small_b}
\end{claim}


{\it Proof.} We will provide a detailed proof of Claim~\ref{small_b} for $b=5$ as the most laborious case and outline similar steps for $b \in \{3,4\}$.
Specifically, for $b=5$, we want to prove that, for all $n \geq 17$:
$$
n^np_{5,n}=(n-5)^n + 5n(n-5)^{n-1} + 25\frac{n(n-1)}{2}(n-5)^{n-2} + 
$$ 
$$
125\frac{n(n-1)(n-2)}{6}(n-5)^{n-3} + 625\frac{n(n-1)(n-2)(n-3)}{24}(n-5)^{n-4} < 
$$
$$
(n-6)^n + 6n(n-6)^{n-1} + 18n(n-1)(n-6)^{n-2} + 36n(n-1)(n-2)(n-6)^{n-3} + 
$$
$$ 
54n(n-1)(n-2)(n-3)(n-6)^{n-4} + \frac{324}{5}n(n-1)(n-2)(n-3)(n-4)(n-6)^{n-5}=p_{6,n}n^n.
$$

For $b\in\{3,4,5\}$, consider the difference $n^n(p_{b+1,n}-p_{b,n})$ and divide it by $(n-b)^n$. Since  
$$
\left(1 - \frac{1}{n - b}\right)^{n-b} > e^{-1}\left(1 - \frac{1}{2(n-b)} - \frac{1}{(n-b)^2}\right),
$$
$p_{6,n}>p_{5,n}$ follows from
$$
\frac{899}{5} - \frac{523e}{8} + \frac{20531\cdot6 - 9025\cdot5e}{60(n - 5)} + \frac{80527\cdot4 - 24625\cdot5e}{40(n - 5)^2} + 
$$ 
$$
\frac{11009\cdot12 - 63125e}{12(n - 5)^3} - \frac{20529 - 3125\cdot5e}{5(n - 5)^4} - \frac{217\,291}{10(n - 5)^5} - \frac{156\,627}{10 (n - 5)^6} - \frac{3125}{(n - 5)^7} > 0,
$$
$p_{5,n} > p_{4, n}$ similarly follows from 
$$\frac{523}{8} - \frac{71e}{3} + \frac{21979 - 168\cdot48e}{48(n-4)} + \frac{21683 - 1120 \cdot 8e}{24(n-4)^2} + \frac{61 - 256 \cdot 48e}{48(n-4)^3} -$$ $$- \frac{18091}{12(n-4)^4} - \frac{14723}{12(n-4)^5} - \frac{256}{(n-4)^6} > 0,$$
and $p_{4,n} > p_{3,n}$ follows from 

$$\frac{142 - 51e}{6} + \frac{511 - 189e}{6(n-3)} - \frac{162e - 217}{6(n-3)^2} - \frac{745}{6(n-3)^3} - \frac{731}{6(n-3)^4} - \frac{27}{(n-3)^5} > 0.$$

To prove these inequalities we consider separately two parts of the summations: the first part contains the constant term and a fraction of the next one, and the second part contains the rest. In particular, when $b=5$, we prove that, for an appropriate choice of $C>0$, the following inequalities hold true:
\begin{equation}
\frac{899}{5} - \frac{523e}{8} + \frac{20531\cdot6 - 9025\cdot5e - C}{60(n - 5)} > 0.
\label{ineq_1}
\end{equation}
\begin{equation}
\begin{aligned}
C + \frac{3(80527\cdot4 - 24625\cdot5e)}{2(n - 5)} + \frac{5(11009\cdot12 - 63125e)}{(n - 5)^2} - \frac{12(20529 - 3125\cdot5e)}{(n - 5)^3} - \\
-\frac{6\cdot217\,291}{(n - 5)^4} - \frac{6\cdot156\,627}{(n - 5)^5} - \frac{60\cdot3125}{(n - 5)^6} > 0,
\end{aligned}
\label{ineq_2}
\end{equation} 
If so, summation of~(\ref{ineq_1}) with~(\ref{ineq_2}) divided by $60(n-5)$ gives the former inequality.

It is easy to see that, for $C := 2300$, both functions on the left sides of~(\ref{ineq_1}),~(\ref{ineq_2}) increase in $n$ and are positive for $n = 20$. So, for $b = 5$ and $n \geq 20$, this finishes the proof. In the remaining three cases, it can be computed that
$$
p_{5,17}17^{17} < 3.387\cdot 10^{20} < 3.389\cdot 10^{20}<p_{6,17}17^{17},
$$
$$
p_{5,18}18^{18}< 1.619\cdot 10^{22} < 1.622\cdot 10^{22} < p_{6,18}18^{18},
$$
$$
p_{5,19}19^{19}<8.176\cdot 10^{23} < 8.199\cdot 10^{23}< p_{6,19}19^{19}.
$$
It remains to check manually that $p_{5,n} > p_{6,n}$ for every $n \leq 16$ (e.g.,  $p_{5,16}16^{16} > 7.505\cdot 10^{18} > 7.503\cdot 10^{18} > p_{6,16}16^{16}$), and we skip these simple computations. For the cases $b=3$ and $b=4$ our method proves that $p_{b,n} < p_{b+1, n}$ for $n \geq 13$ and $n \geq 17$ (as opposed to claimed $n \geq 11$ and $n \geq 14$) respectively. The proof for the rest values of $n$ can be obtained by the same numerical approach. For $b=1$ and $b=2$, the statement can be easily checked. $\quad\Box$

\begin{claim}
Let $b\geq 6$ and $b\in\{n-5,n-4,n-3,n-2,n-1\}$. Then $p_{b,n} > p_{b+1, n}$.
\label{p_bn_for_b_equals_n_minus_5}
\end{claim}

{\it Proof.} First, $p_{n-1,n}>p_{n,n}$ is obvious since $p_{n,n}=0$ while $p_{n-1,n}$ is positive.

Second, 
$$
 n^n(1-p_{n-2,n})={n\choose 2}4(n-2)^{n-2}+2n(n-2)^{n-1}+(n-2)^n,
$$
$$
 n^n(1-p_{n-1,n})=n(n-1)^{n-1}+(n-1)^n.
$$
Therefore, $p_{n-2,n}>p_{n-1,n}$ if and only if
$$
 2\left(1+\frac{2}{n-2}\right)\left(1+\frac{1}{n-2}\right)+2\left(1+\frac{2}{n-2}\right)+1<
$$
$$
 \left[\left(1+\frac{2}{n-2}\right)\left(1+\frac{1}{n-2}\right)+\left(1+\frac{1}{n-2}\right)^2\right]\left(1+\frac{1}{n-2}\right)^{n-2}.
$$
Denote $x:=\frac{1}{n-2}$. Since $n\geq 8$, we get that $\left(1+\frac{1}{n-2}\right)^{n-2}>\frac{5}{2}$. Therefore, the desired inequality follows from 
$$
 \frac{5}{2}\left[\left(1+2x\right)\left(1+x\right)+\left(1+x\right)^2\right]=\frac{5}{2}(3x^2+5x+2)>4x^2+10x+5=2(1+2x)(1+x)+2(1+2x)+1.
$$

It remains to consider $b\in\{n-5,n-4,n-3\}$.

From~(\ref{p_above_n_over_2}), we get
$$
n^n(p_{b+1,n}-p_{b,n})=
$$
$$
\frac{n(n-1)\ldots(b+1)}{(n-b-1)!}\left((n-b)^{n-b-1}b^b -
(n-b-1)^{n-b-1}(b+1)^b\right)+n^n(p_{n-b,n}-p_{n-b-1,n}).
$$
Let us divide both parts of the latter inequality by $(b+1)^n$, substitute 
$$
\sum_{i=0}^{n-b-1}{n\choose i}\left(\frac{n-b}{n}\right)^i\left(1-\frac{n-b}{n}\right)^{n-i}\text{ and }\,\,\sum_{i=0}^{n-b-2}{n\choose i}\left(\frac{n-b-1}{n}\right)^i\left(1-\frac{n-b-1}{n}\right)^{n-i}
$$ 
in $p_{n-b,n}$ and $p_{n-b-1,n}$ respectively and apply the inequality 
$$
\left(1 - \frac{1}{b+1}\right)^k=
\left(1 - \frac{1}{b+1}\right)^{b+1}\left(1 - \frac{1}{b+1}\right)^{k-b-1}<
 e^{-1}\left(1-\frac{k-b-1}{b+1}\right)
$$ 
which is true for all $6\leq b\leq k \leq n$. 

Then, for $b=n-5$, we get


$$
\frac{n^n}{(n-4)^n}(p_{n-4,n}-p_{n-5,n})<\frac{1097}{12e}-\frac{103}{3} +\frac{18649/(24e)-824/3}{n-4}+
$$
$$
\frac{4705/(2e)-2240/3}{(n-4)^2}+\frac{832/3-5225/(8e)}{(n-4)^3}+\frac{24625/(12e)-256}{(n-4)^4}+\frac{625/e}{(n-4)^5}.
$$
The bound is negative for $n=28$. Moreover, all the terms except the constant are positive for $n > 4$, which implies the negativeness of the bound for all $n \geq 28$.

Let $n=11$ and $b=6$. Then:
$$p_{6,11}11^{11} = {11 \choose 0} 6^0 5^{11} + {11 \choose 1} 6^1 5^{10} + {11 \choose 2} 6^2 5^9 + {11 \choose 3} 6^3 5^8 + {11 \choose 4} 6^4 5^7 + {11 \choose 5} 6^5 5^6 = $$
$$ = 48828125+
644531250+
3867187500+ $$ $$ +
13921875000+
33412500000+
56133000000 > 1.08 \times 10^{11} $$
$$p_{7,11} 11^{11} = {11 \choose 0} 7^0 4^{11} + {11 \choose 1} 7^1 4^{10} + {11 \choose 2} 7^2 4^9 + {11 \choose 3} 7^3 4^8 + {11 \choose 4} 7^4 4^7 + $$ $$ + {11 \choose 5} 7^5 4^6 + {11 \choose 6} 7^6 4^5
= 4194304+
80740352+
706478080+$$ $$ +
3709009920+
12981534720+
31804760064+
55658330112 < 1.05 \times 10^{11} ,
$$
therefore, $p_{b, n} > p_{b + 1, n}$ for $n = 11$ and $b = n - 5 = 6$.

Now let $n = 12$ and $b = 7$.
$$p_{7,12} 12^{12} = {12 \choose 0} 7^0 5^{12} + {12 \choose 1} 7^1 5^{11} + {12 \choose 2} 7^2 5^{10} +  {12 \choose 3} 7^3 5^9 + {12 \choose 4} 7^4 5^8 + $$
$$ + {12 \choose 5} 7^5 5^7  + {12 \choose 6} 7^6 5^6 =244140625+
4101562500+
31582031250+
147382812500+ $$ $$ +
464255859375+ 
1039933125000+
1698557437500 > 3.3 \times 10^{12} 
$$
$$p_{8,12} 12^{12} = {12 \choose 0} 2^0 + {12 \choose 1} 2^1 + {12 \choose 2} 2^2 +  {12 \choose 3} 2^3 + {12 \choose 4} 2^4 + {12 \choose 5} 2^5 + $$
$$ + {12 \choose 6} 2^6 + {12 \choose 7} 2^7 = 16777216+
402653184+
4429185024+
29527900160+ $$ $$ +
132875550720+
425201762304+ 
992137445376+
1700807049216 < 3.3 \times 10^{12}
,$$
therefore, $p_{b, n} > p_{b+1, n}$ for $n=12$ and $b = n-5 = 7$. The calculations for the remaining cases ($13 \leq n \leq 27$) follow the exact same pattern (see \url{https://github.com/daniildmitriev/on_monotonicity_ramanujan/blob/master/verification_of_claim_8.txt}).

In the same way, for $b=n-4$, we get 
$$
\frac{n^n}{(n-3)^n}(p_{n-3,n}-p_{n-4,n})<\frac{103}{3e}-13 +\frac{515/(3e)-117/2}{n-3}+
$$
$$
\frac{848/(3e)-153/2}{(n-3)^2}+\frac{592/(3e)-27}{(n-3)^3}+\frac{64/e}{(n-3)^4}.
$$
The bound is negative for $n=21$, and all the terms except the constant are positive for $n > 3$, which implies the negativeness of the bound for all $n \geq 21$. The calculations for the remaining cases ($10 \leq n \leq 20$) are in the aforementioned file.

For $b = n-3$, we get
$$
\frac{n^n}{(n-2)^n}(p_{n-2,n}-p_{n-3,n})<\frac{13}{e}-5 +\frac{65/(2e)-10}{n-2}+
\frac{51/(2e)-4}{(n-2)^2}+\frac{9/e}{(n-2)^4}.
$$
The bound is negative for $n=14$, and all the terms except the constant are positive for $n > 2$, which implies the negativeness of the bound for all $n \geq 14$. The calculations for the remaining cases ($9 \leq n \leq 13$) are in the aforementioned file.


$\quad\Box$


\section*{Appendix C. Proof of Claim~\ref{th1_main_claim}}
\subsection*{Proof of inequality~(\ref{small_b_ineq})}

Let $6 \leq b \leq \frac{n-2}{3}$. Clearly,
$$
 \frac{bn(\partial^{4} g/\partial z^{4})\left(1-\frac{b}{n}\right)}{5\left(\frac{b+1}{n}\right)^{b-4}\left(1-\frac{b+1}{n}\right)^{n-b-2}}
 =
 \frac{bn^4}{5(b+1)(n-b-1)^2}\left(\frac{b}{b+1}\right)^{b-5}\left(\frac{n-b}{n-b-1}\right)^{n-b-4}\times
$$
$$
\left[3(n-b)^2\left(\frac{b}{n}\right)^2-2(n-b)\left(\frac{b}{n}\right)\left(23\left(1-\frac{b}{n}\right)^2-13\left(1-\frac{b}{n}\right)+3\right)+24\left(1-\frac{b}{n}\right)^4\right]. 
$$
As $ \ln(1+x) \leq x$ for all $x>-1$, and $e^x < 1 + 2x$ for all $x\in(0,1)$, we get
$$
\left(\frac{b}{b+1}\right)^{b-5}\left(\frac{n-b}{n-b-1}\right)^{n-b-4}=\exp\left((b-5)\ln\left(1-\frac{1}{b+1}\right)+(n-b-4)\ln\left(1+\frac{1}{n-b-1}\right)\right)
$$
$$
\leq\exp\left(-\frac{b-5}{b+1}+\frac{n-b-4}{n-b-1}\right)=\exp\left(\frac{6}{b+1}-\frac{3}{n-b-1}\right)<\exp\left(\frac{6}{b+1}\right)<1+\frac{12}{b+1}.
$$
Then,
$$
\frac{bn(\partial^{4} g/\partial z^{4})\left(1-\frac{b}{n}\right)}{5\left(\frac{b+1}{n}\right)^{b-4}\left(1-\frac{b+1}{n}\right)^{n-b-2}}<\frac{bn^4(b+13)}{5(b+1)^2(n-b-1)^2}\left(3(n-b)^2\left(\frac{b}{n}\right)^2+24\left(1-\frac{b}{n}\right)^4\right).
$$
It remains to prove that the expression to the right in the last inequality is smaller than $P_{b,n}$, which, in turn, follows from
\begin{equation}
5(b+1)^2(n-b-1)^2P_{b,n} > 3bn^2(b+13)(b^2+8)(b-n)^2.
\label{ineq_22}
\end{equation}

Set
$$
 R_{b,n}=156b^2+12bn^2-105bn+94b+24n^2-48n+24,
$$
and so
$$
 P_{b,n}=12b^5-16b^4n+64b^4+4b^3n^2-71b^3n+138b^3+16b^2n^2-112b^2n+R_{b,n}.
$$
Surely, for $6\leq \frac{n-2}{3}$, $3n^2-48n$ is positive. Moreover, $16(3b-n)^2+5(b-n)^2=149b^2-106bn+21n^2$. Therefore, $R_{b,n}>0$.

From this,
\begin{equation}
P_{b,n} > b^2(12b^3-16b^2n+64b^2+4bn^2-71bn+138b+16n^2-112n).
\label{P_bn_small_b}
\end{equation}

First, assume that $b\geq 39$.

Let us prove that $\frac{n-2}{3}\geq b$ implies $(b+1)^2(n-b-1)^2>(b^2+8)(b-n)^2$. Indeed, $(b+1)^2(n-b-1)^2-(b^2+8)(b-n)^2=(n-b)^2(2b-7)-2(n-b)(b+1)^2+(b+1)^2$. The function $x^2(2b-7)-2x(b+1)^2+(b+1)^2$ increases when $x>\frac{(b+1)^2}{2b-7}$ and is positive in the point $x=2b+2>\frac{(b+1)^2}{2b-7}$. Therefore, $(n-b)^2(2b-7)-2(n-b)(b+1)^2+(b+1)^2>0$ for all $n\geq 3b+2$ as desired.

Now, it remains to prove that
$$
5b(12b^3-16b^2n+64b^2+4bn^2-71bn+138b+16n^2-112n) - 3n^2(b+13) > 0.
$$
For this polynomial of $n$ to be positive for all (real) $n\geq 3b+2$, its roots must be less than $3b+2$. This leads to the inequality $\frac{-B + \sqrt{B^2 - 4AC}}{2A} < 3b+2$, where $A = 20b^2+77b-39$, $B = -80b^3-355b^2-560b$,
$C = 60b^4+320b^3+690b^2$. The latter inequality is equivalent to
$$
4(20b^2+77b-39)(28b^3-1047b^2-1280b-156) > 0,
$$
which is true for $b \geq 39$.

Now, let us prove the inequality~(\ref{ineq_22}) for $6\leq b \leq 38$, $n \geq 158$. Since $b^2+8<(b+1)^2$, it is enough to show that
\begin{equation}
P_{b,n} > \frac{3}{5}b(b+13)n^2\left(1 + \frac{1}{n-b-1}\right)^2.
\label{P_bn_small_b_second}
\end{equation}
Since, in the considered range of the parameters,
$$
 \frac{3}{5}\frac{b+13}{b}\left(1 + \frac{1}{n-b-1}\right)^2<2,
$$
the inequality~(\ref{P_bn_small_b_second}) follows from~(\ref{P_bn_small_b}) and
$$
12b^3-16b^2n+64b^2+4bn^2-71bn+138b+16n^2-112n > 2n^2
$$
(the latter is straightforward since $n\geq 4b+6$).

In all the remaining cases ($6 \leq b \leq 38$ and $3b+2 \leq n \leq 157$), we prove inequality~(\ref{small_b_ineq}) by verifying that $ {\sf LHS} := \frac{bn(\partial^{4} g/\partial z^{4})\left(1-\frac{b}{n}\right)}{5\left(\frac{b+1}{n}\right)^{b-4}\left(1-\frac{b+1}{n}\right)^{n-b-2}} < P_{b,n} =: {\sf RHS}$. We compute the values of both ${\sf LHS}$ and ${\sf RHS}$ manually and obtain the following inequalities.

\begin{enumerate}
\item For $b = 6$ and $n = 3b + 2 = 14$:
${\sf LHS} < 4452.57 < 11028 = {\sf RHS}.$
\item For $b = 6$ and $n = 3b + 3 = 15$:
${\sf LHS} < 4560.47 < 33222 = {\sf RHS}.$
\item For $b = 6$ and $n = 156$:
${\sf LHS} < 0 < 31230372 = {\sf RHS}.$
\item For $b = 6$ and $n = 157$:
${\sf LHS} < 0 < 31670358 = {\sf RHS}.$
\item For all the remaining cases, please see \url{https://github.com/daniildmitriev/on_monotonicity_ramanujan/blob/master/verification_of_inequality_11.txt}
\end{enumerate}

\subsection*{Proof of inequality~(\ref{negativeness})}
Here, $\frac{n-1}{3} \leq b \leq \frac{n}{2}$.
Since 
$$
\frac{\partial^{4} g}{\partial z^{4}}\left(1-\frac{b+1}{n}\right)=
 \left(\frac{b+1}{n}\right)^{b-5}\left(1-\frac{b+1}{n}\right)^{n-b-4}
 \biggl[3(n-b-1)^2\left(\frac{b+1}{n}\right)^2-
$$
$$
 2(n-b-1)\left(\frac{b+1}{n}\right)\left(23\left(1-\frac{b+1}{n}\right)^2+7\left(1-\frac{b+1}{n}\right)-1\right)+
$$
$$
96\left(1-\frac{b+1}{n}\right)^3+24\left(1-\frac{b+1}{n}\right)^4\biggr]>
\left(\frac{b+1}{n}\right)^{b-5}\left(1-\frac{b+1}{n}\right)^{n-b-4}\times
$$
$$\biggl[3(n-b-1)^2\left(\frac{b+1}{n}\right)^2-2(n-b-1)\left(\frac{b+1}{n}\right)\left(23\left(1-\frac{b+1}{n}\right)^2+7\left(1-\frac{b+1}{n}\right)\right)\biggr]=
$$
\begin{equation}
\left(\frac{b+1}{n}\right)^{b-4}\left(1-\frac{b+1}{n}\right)^{n-b-2}(3bn+46b-57n+46),
\label{4_from_below}
\end{equation}
 the negativeness of~(\ref{negativeness}) will follow from
\begin{equation}
5P_{b,n} < bn(3bn+46b-57n).
\label{medium_b_ineq}
\end{equation}
We want to show that $5P_{b,n} - bn(3bn+46b-57n) < 0$, or, as in the proof of inequality~(\ref{small_b_ineq}), that the segment $[2b, 3b + 1]$ lies inside $(n_1, n_2)$, where $n_1$ and $n_2$ are the smallest and the largest roots of this polynomial of $n$ respectively. Since $n_2 > \frac{-B+\sqrt{B^2-4AC}}{2A}$, where $A=20b^2+77b+129$, $B=-80b^3-355b^2-606b$, $C=60b^4+320b^3+690b^2+830b$, the inequality $n_2 > 3b+1$ follows from
$$
4(20b^2+77b+129)(12b^3-160b^2-1075b-129) > 0,
$$
which is true for $b \geq 19$. Similarly it can be proven that $n_1 < 2b$ for $b \geq 19$.

For the remaining cases ($6 \leq b \leq 18$, $2b \leq n \leq 3b+1$) we verify that $ {\sf LHS} := 5 P_{b,n} < \frac{bn(\partial^{4} g/\partial z^{4})\left(1-\frac{b}{n}\right)}{5\left(\frac{b+1}{n}\right)^{b-4}\left(1-\frac{b+1}{n}\right)^{n-b-2}} =: {\sf RHS}$ manually.
\begin{enumerate}
\item For $b = 6$ and $n = 2b = 12$:
${\sf LHS} = -279660 < 907 < {\sf RHS}.$
\item For $b = 6$ and $n = 2b+1= 13$:
${\sf LHS} = -291570 < 919 < {\sf RHS}.$
\item For $b = 6$ and $n = 2b+2= 14$:
${\sf LHS} = -288120 < 915 < {\sf RHS}.$
\item For $b = 6$ and $n = 2b+3= 15$:
${\sf LHS} = -269310 < 900 < {\sf RHS}.$
\item For $b = 6$ and $n = 2b+4 = 16$:
${\sf LHS} = -235140 < 876 < {\sf RHS}.$
\item For $b = 6$ and $n = 2b+5 = 17$:
${\sf LHS} = -185610 < 847 < {\sf RHS}.$
\item For $b = 6$ and $n = 2b+6 = 18$:
${\sf LHS} = -120720 < 814 < {\sf RHS}.$
\item For $b = 6$ and $n = 2b+7 = 3b+1 = 19$:
${\sf LHS} = -40470 < 779 < {\sf RHS}.$
\item For all the remaining cases, please see \url{https://github.com/daniildmitriev/on_monotonicity_ramanujan/blob/master/verification_of_inequality_13.txt}
\end{enumerate}

\subsection*{Proof of inequality~(\ref{negativeness2})}
Let $5 \leq b < \frac{n}{3}$. Since $\ln(1-x)>-x-x^2/2-x^3/2$ for $x\in(0,1/2]$, and $\ln(1+x)>x-x^2/2$ for positive $x$, we get
$$
\frac{n-b}{n-b-1} - \left(\frac{b}{b+1}\right)^{b}\left(\frac{n-b}{n-b-1}\right)^{n-b}=
$$
$$
1+\frac{1}{n-b-1} -\mathrm{exp}\left(b\ln\left(1-\frac{1}{b+1}\right)+(n-b)\ln\left(1+\frac{1}{n-b-1}\right)\right)<
$$
$$
1+\frac{1}{n-b-1} -\mathrm{exp}\left(b\left(-\frac{1}{b+1}-\frac{1}{2(b+1)^2}-\frac{1}{2(b+1)^3}\right)+(n-b) \left(\frac{1}{n-b-1}-\frac{1}{2(n-b-1)^2}\right)\right)
$$
$$
=1+\frac{1}{n-b-1} -\mathrm{exp}\left(\frac{1}{2(b+1)}-\frac{1}{2(b+1)^3}+\frac{1}{2(n-b-1)}-\frac{1}{2(n-b-1)^2}\right)<
$$
$$
=1+\frac{1}{n-b-1} -\mathrm{exp}\left(\frac{1}{2(b+1)}+\frac{1}{2(n-b-1)}-\frac{1}{2(n-b-1)^2}\right)<
$$
$$
=1+\frac{1}{n-b-1} -\left(1+\frac{1}{2(b+1)}+\frac{1}{2(n-b-1)}-\frac{1}{2(n-b-1)^2}\right)=
$$
$$
< - \frac{1}{2(b+1)} + \frac{n-b}{2(n-b-1)^2}.
$$
From (\ref{4_from_below}), we get that (\ref{negativeness2}) follows from
$$
5P_{b,n}-bn(3bn+46b-57n)+60(b+1)^3(3nb-2b^2+3n-3b-n^2-1) < 0.
$$
We should prove that $An^2+Bn+C<0$ where
\begin{align*}
A=-40(b+1)^3+27(b+1)^2+3(b+1)+70,\\
B=100(b+1)^4-35(b+1)^3-21(b+1)^2-78(b+1)-26,\\
C=-60(b+1)^5+80(b+1)^4+10(b+1)^3+30(b+1)^2.
\end{align*}
Since, for considered values of $b$, we have $A < 0$, the inequality is true if the biggest root of the polynomial is less than $3b$. It is enough to check that at the point $n = 3b$ the polynomial and its derivative are negative. This brings us to the two inequalities, $9Ab^2+3Bb+C<0$ and $6Ab + B < 0$. Let us verify these inequalities using the fact that $b \geq 5$:
$$9Ab^2 + 3Bb + C = -120 b^5 + 38b^4 + 585b^3 + 1005b^2 - 70b + 60 < b^2 (-120b^3 + 38b^2 + 585b + 1005) < $$
$$ < b^3(-110b^2 + 38b + 585) < b^4(-80b + 38) < 0,$$
$$6Ab + B = -140b^4 - 193b^3 + 96b^2 + 535b - 60 < b^2 (-140b^2 + 96) < 0.$$
This concludes the proof.

\section*{Appendix D. Plots of $p_{b,n}$, $z_{b,n}$ and $\tilde p_{c,b,n}$}
\begin{figure}[h]
\caption{Plot of $p_{b,n}$ for $n = 1000$, for all values of $b$ (left) and for $b \in [n / 3 - 100, n / 3 + 100]$ (right). Red line represents $n = 3b + 1$ -- point of monotonicity change.}
\centering
\includegraphics[width=0.99\textwidth]{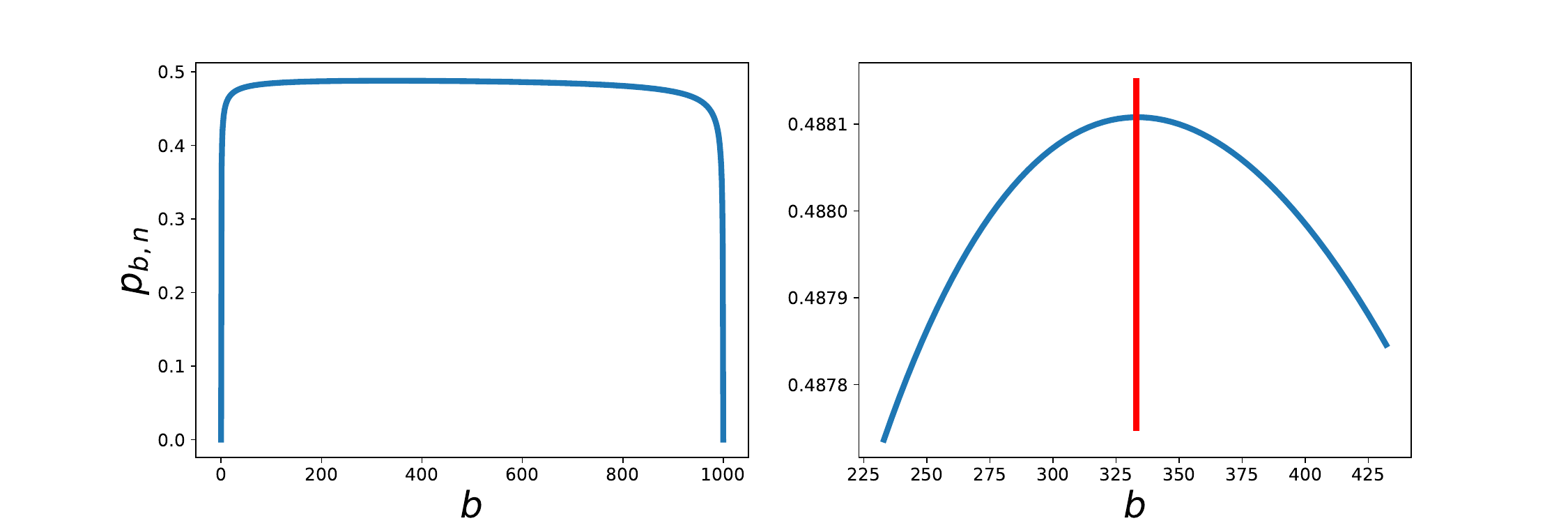}
\end{figure}
\begin{figure}[h]
\caption{Plot of $z_{b,n}$ for $n = 1000$, for all values of $b$ (left) and for $b \leq \sqrt{n}$ (right).}
\centering
\includegraphics[width=0.99\textwidth]{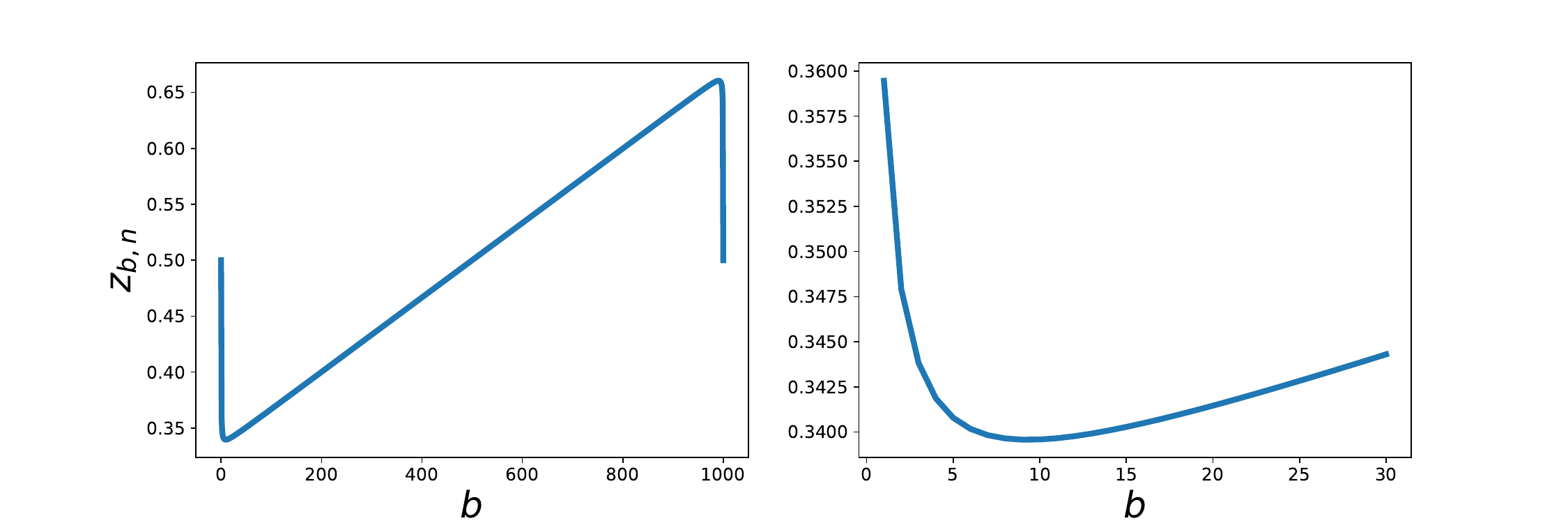}
\end{figure}
\begin{figure}[h]
\caption{Plot of $p_{c,b,n}$ for $n = 1000$, for $c \in \{0, \frac{1}{2}, \frac{2}{3}, 1\}$ for $0 \leq b \leq n$ (left) and for $\frac{n}{10} \leq b \leq \frac{9n}{10}$ (right).}
\centering
\includegraphics[width=0.99\textwidth]{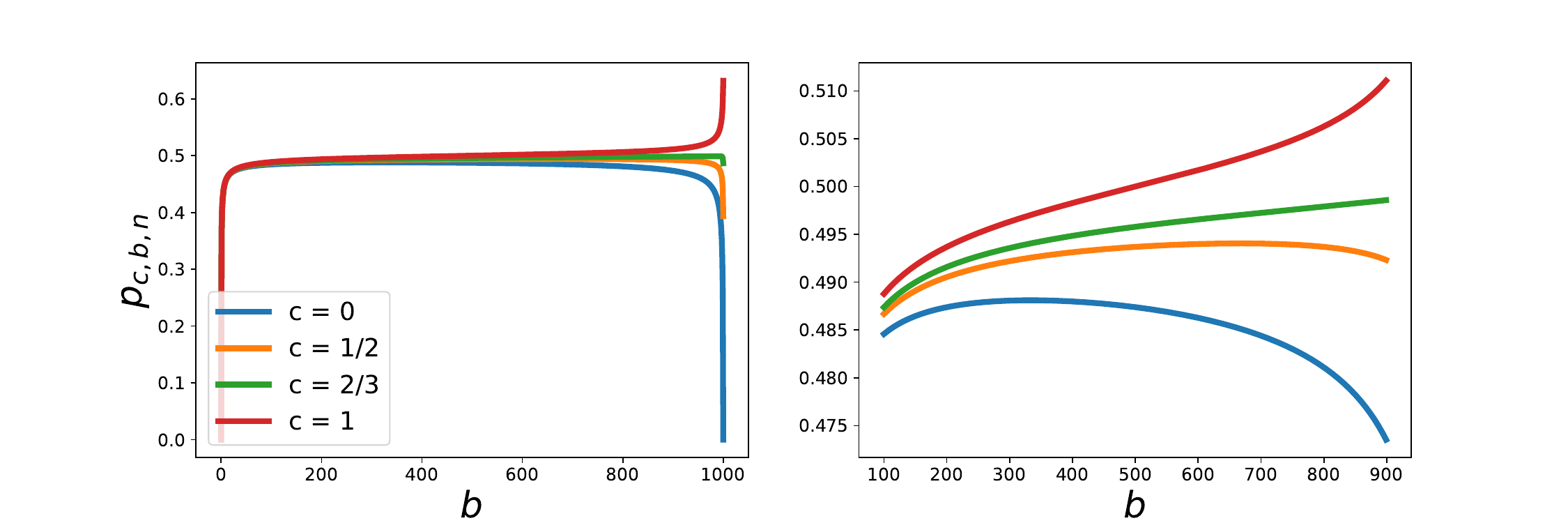}
\end{figure}
\end{document}